\documentclass[reqno,12pt]{amsart}
\usepackage{amsmath}

\usepackage{amsfonts}

\usepackage{amssymb}

\usepackage{eufrak}

\usepackage{amscd}

\usepackage{amsthm}

\usepackage{amstext}

\usepackage[all]{xy}

   \newcommand{\Hom}{\operatorname{Hom}}
   
 \newcommand{\Ext}{\operatorname{Ext}}
\newcommand{\Ad}{\operatorname{Ad}}

\newcommand{\id}{\operatorname{id}}

\newcommand{\im}{\operatorname{im}}

   \theoremstyle{plain}
   \newtheorem{thm}{Theorem}
   
   \newtheorem{lemma}[thm]{Lemma}
   
   \theoremstyle{definition}

   \theoremstyle{remark}

\title[Extensions]{The group of unital  $C^*$-extensions}
\author{Vladimir Manuilov}
\address{Department of Mechanics and Mathematics, Moscow State
University, Leninskie Gory, Moscow, 119992, Russia}
\email{manuilov@mech.math.msu.su}

\author{Klaus Thomsen}
\address{IMF, Department of Mathematics, Ny Munkegade, 8000 Aarhus C,
Denmark}
\email{matkt@imf.au.dk}

\begin{document}




\begin{abstract} Let $A$ and $B$ be separable $C^*$-algebras, $A$ unital
and $B$ stable. It is shown that there is a natural six-terms exact
sequence which relates the group which arises by considering \emph{all}
semi-split extensions of $A$ by $B$ to the group which arises by
restricting the attention to \emph{unital} semi-split extensions of $A$ by
$B$. The six-terms exact sequence is an unpublished result of G.
Skandalis.

\end{abstract}

\maketitle




Let $A,B$ be separable $C^*$-algebras, $B$ stable. As is well-known the
$C^*$-algebra extensions of $A$ by $B$ can be identified with $\Hom
(A.Q(B))$, the set of $*$-homomorphisms $A \to Q(B)$ where $Q(B) = M(B)/B$
is the generalized Calkin algebra. Two extensions $\varphi, \psi : A \to
Q(B)$ are \emph{unitarily equivalent} when there is a unitary $u \in M(B)$
such that $\Ad q(u) \circ \psi = \varphi$, where $q : M(B) \to Q(B)$ is
the quotient map. The unitary equivalence classes of extensions of $A$ by
$B$ have an abelian semi-group structure thanks to the stability of $B$:
Choose isometries $V_1,V_2 \in M(B)$ such that $V_1V_1^* + V_2V_2^* = 1$,
and define the sum $\varphi \oplus \psi : A \to Q(B)$ of $\varphi, \psi
\in\Hom (A,Q(B))$ by
\begin{equation}\label{B1}
(\psi \oplus \varphi)(a) = \Ad q(V_1) \circ \psi(a)     + \Ad q(V_2)
\circ \varphi(a) .
\end{equation}
The isometries, $V_1$ and $V_2$, are fixed in the following. An extension
$\varphi : A \to Q(B)$ is \emph{split} when there is a $*$-homomorphisms
$\pi : A \to M(B)$ such that $\varphi = q \circ \pi$. To trivialize the
split extensions we declare two extensions  $\varphi, \psi : A \to Q(B)$
to be \emph{stably equivalent} when there there is a split extension $\pi$
such that $\psi \oplus \pi$ and $\varphi \oplus \pi$ are unitarily
equivalent. This is an equivalence relation because the sum (\ref{B1}) of
two split extensions is itself split. We denote by $\Ext(A,B)$ the
semigroup of stable equivalence classes of extensions of $A$ by $B$. It
was proved in \cite{Th}, as a generalization of results of Kasparov, that
there exists an \emph{absorbing split extension} $\pi_0 : A \to Q(B)$,
i.e. a split extension with the property that $\pi_0 \oplus \pi$ is
unitarily equivalent to $\pi_0$ for every split extension $\pi$. Thus two
extensions $\varphi,\psi$ are stably equivalent if and only if $\varphi
\oplus \pi_0$ and $\psi \oplus \pi_0$ are unitarily equivalent. The
classes of stably equivalent extensions of $A$ by $B$ is an abelian
semigroup $\Ext(A,B)$ in which any split extension (like $0$) represents
the neutral element. As is well-documented the semi-group is generally not
a group, and we denote by
 $$
\Ext^{-1}(A,B)
 $$
the abelian group of invertible elements in $\Ext(A,B)$. It is also
well-known that this group is one way of describing the $KK$-groups of
Kasparov. Specifically, $\Ext^{-1}(A,B) = KK(SA,B) = KK (A,SB)$.

\smallskip

Assume now that $A$ is unital. It is then possible, and sometimes even
advantageous, to restrict attention to unital extensions of $A$ by $B$,
i.e. to short exact sequences
\begin{equation}\label{B2}
\begin{xymatrix}{
0 \ar[r]  & B \ar[r] & E \ar[r]  &  A \ar[r]  & 0}
\end{xymatrix}
\end{equation}
of $C^*$-algebras with $E$ is unital, or equivalently to
$*$-homomorphisms $A \to Q(B)$ that are unital. The preceding definitions
are all amenable to such a restriction, if done consistently.
Specifically, we say that a unital extension $\varphi : A \to Q(B)$ is
\emph{unitally split} when there is a unital $*$-homomorphism $\pi : A \to
M(B)$ such that $\varphi = q \circ \pi$. The sum (\ref{B1}) of two unital
extensions is again unital, and we say that two unital extensions
$\varphi, \psi : A \to Q(B)$ are \emph{unitally stably equivalent} when
there there is a unital split extension $\pi$ such that $\psi \oplus \pi$
and $\varphi \oplus \pi$ are unitarily equivalent. It was proved in
\cite{Th} that there always exists a \emph{unitally absorbing split
extension} $\pi_0 : A \to Q(B)$, i.e. a unitally split extension with the
property that $\pi_0 \oplus \pi$ is unitarily equivalent to $\pi_0$ for
every unitally split extension $\pi$. Thus two unital extensions
$\varphi,\psi$ are unitally stably equivalent if and only if $\varphi
\oplus \pi_0$ and $\psi \oplus \pi_0$ are unitarily equivalent. The
classes of unitally stably equivalent extensions of $A$ by $B$ is an
abelian semi-group which we denote by $\Ext_{\text{unital}}(A,B)$. The
unitally absorbing split extension $\pi_0$, or any other unitally split
extension, represents the neutral element of $\Ext_{\text{unital}}(A,B)$,
and we denote by
 $$
\Ext_{\text{unital}}^{-1}(A,B)
 $$
the abelian group of invertible elements in $\Ext_{\text{unital}}(A,B)$.
As we shall see there is a difference between
$\Ext_{\text{unital}}^{-1}(A,B)$ and $\Ext^{-1}(A,B)$ arising from the
fact that while the class in $\Ext^{-1}(A,B)$ of a unital extension $A \to
Q(B)$ can not be changed by conjugating it with a unitary from $Q(B)$, its
class in $\Ext_{\text{unital}}^{-1}(A,B)$ can. In a sense the main result
of this note is that this is the only way in which the two groups differ.

Note that there is a group homomorphism
 $$
\Ext_{\text{unital}}^{-1}(A,B) \to \Ext^{-1}(A,B),
 $$
obtained by forgetting the word 'unital'. It will be shown that this
forgetful map fits into a six-terms exact sequence
\begin{equation}\label{B3}
\begin{xymatrix}{
K_0(B) \ar[r]^-{u_0} & \Ext^{-1}_{\text{unital}}(A,B) \ar[r] &
\Ext^{-1}(A,B) \ar[d]^{i_1^*} \\ \Ext^{-1}(A,SB) \ar[u]^-{i_0^*}  &
\Ext^{-1}_{\text{unital}}(A,SB) \ar[l] & K_1(B) \ar[l]^-{u_1} }
\end{xymatrix}
\end{equation}
where $SB$ is the suspension of $B$, i.e. $SB = C_0(0,1) \otimes B$, and
the maps $u_k$ and $i_k^*, k = 0,1$, will be defined shortly. This
six-terms exact sequence is mentioned in 10.11 of \cite{S}, but the proof
was never published.

\smallskip

Fix a unitally absorbing $*$-homomorphism $\alpha_0 : A \to M(B)$, which
exists by Theorem 2.4 of \cite{Th}. It follows then from Theorem 2.1 of
\cite{Th} that $\alpha = q \circ \alpha_0$ is a unitally absorbing split
extension as defined above.

\begin{lemma}\label{B31} The $*$-homomorphisms $\Ad V_1 \circ \alpha_0 :
A \to M(B)$ and $\left( \begin{smallmatrix} \alpha_0  & \\ &
0\end{smallmatrix} \right) : A \to M\left(M_2(B)\right)$ are both
absorbing.
\begin{proof} There is a $*$-isomorphism $ M\left(M_2(B)\right) =
M_2\left(M(B)\right) \to M(B)$ given by
 $$
\left( \begin{smallmatrix}m_{11}  & m_{12} \\  m_{21}& m_{22}
\end{smallmatrix} \right) \mapsto V_1m_{11}V_1^* + V_1 m_{12}V_2^* +
V_2m_{21}V_1^* + V_2m_{22}V_2^*,
 $$
which sends $M_2(B)$ to $B$ and $\left( \begin{smallmatrix} \alpha_0  &
\\ & 0\end{smallmatrix} \right)$ to $\Ad V_1 \circ \alpha_0$, so it
suffices to show that the latter is an absorbing $*$-homomorphism. By
definition, cf. Definition 2.6 of \cite{Th}, we must show that the unital
$*$-homomorphism $A \oplus \mathbb C \ni (a,\lambda) \mapsto
V_1\alpha_0(a)V_1^* + \lambda V_2V_2^*$ is unitally absorbing. For this we
check that it has property 1) of Theorem 2.1 of \cite{Th}. So let $\varphi
: A \oplus \mathbb C \to B$ be a completely positive contraction. Since
$\alpha_0$ has property 1), there is a sequence $\{W_n\}$ in $M(B)$ such
that $\lim_{n \to \infty} W_n^*b = 0$ for all $b \in B$ and $\lim_{n \to
\infty} W_n^* \alpha_0(a)W_n = \varphi(a)$ for all $a \in A$. Since $B$ is
stable there is a sequence $\{S_n\}$ of isometries in $M(B)$ such that
$\lim_{n \to \infty} S_n^*b = 0$ for all $b \in B$. Set
 $$
T_n = V_1W_n + V_2S_n \varphi(0,1)^{\frac{1}{2}}.
 $$
Then $\lim_{n \to \infty} T_nb = 0$ for all $b \in B$, and
 $$
T_n^* \left( V_1\alpha_0(a)V_1^* + \lambda V_2V_2^*\right)T_n =
W_n^*\alpha_0(a)W_n +  \varphi(0,\lambda)
 $$
for all $n$. Since the last expression converges to $\varphi(a,\lambda)$
as $n$ tends to infinity, the proof is complete.
\end{proof}
\end{lemma}

Set
 $$
C_{\alpha} = \left\{ m \in M_2(M(B)): \ m \left( \begin{smallmatrix}
\alpha_0(a)  & \\ & 0\end{smallmatrix} \right) -\left( \begin{smallmatrix}
\alpha_0(a)  & \\ & 0\end{smallmatrix} \right)m \in M_2(B) \ \forall a \in
A \right\}
 $$
and
 $$
A_{\alpha} = \left\{ m \in C_{\alpha}: \ m \left( \begin{smallmatrix}
\alpha_0(a)  & \\ & 0\end{smallmatrix} \right)\in M_2(B) \ \forall a \in A
\right\} .
 $$
We can define a $*$-homomorphism $C_{\alpha} \to \alpha(A)'\cap Q(B)$
such that
 $$
\left( \begin{smallmatrix}m_{11}  & m_{12} \\  m_{21}& m_{22}
\end{smallmatrix} \right) \mapsto q\left(m_{11}\right).
 $$
Then kernel is then $A_{\alpha}$, so we have a $*$-isomorphism
$C_{\alpha}/A_{\alpha} \simeq \alpha(A)'\cap Q(B)$. By Lemma \ref{B31} $
\left( \begin{smallmatrix} \alpha_0  & \\ & 0\end{smallmatrix} \right)$ is
an absorbing $*$-homomorphism, so we conclude from Theorem 3.2 of
\cite{Th} that there is an isomorphism
\begin{equation}\label{B32}
K_1\left(\alpha(A)'\cap Q(B)\right) \simeq KK(A,B) .
\end{equation}
Since the unital $*$-homorphism $\mathbb C \to M(B)$ is unitally
absorbing, this gives us also the well-known isomorphism
\begin{equation}\label{B33}
K_1\left( Q(B)\right) \simeq KK(\mathbb C,B) .
\end{equation}
Let $i: \mathbb C \to A$ be the unital $*$-homorphism. For convenience we
denote the map $K_1\left(\alpha(A)'\cap Q(B)\right) \to K_1(Q(B))$ induced
by the inclusion $\alpha(A)'\cap Q(B) \subseteq Q(B)$ by $i^*$. It is then
easy to check that the isomorphisms (\ref{B32}) and (\ref{B33}) match up
to make the diagram
\begin{equation}\label{B34}
\begin{xymatrix}{
K_1\left(\alpha(A)'\cap Q(B)\right) \ar[d] \ar[r]^-{i^*} & K_1(Q(B))
\ar[d] \\ KK(A,B) \ar[r]^-{i^*} & KK(\mathbb C,B) }
\end{xymatrix}
\end{equation}
commute.

Let $v$ be a unitary in $M_n(Q(B))$. By composing the $*$-homomorphism
$\Ad v \circ \left(1_n \otimes \alpha\right) : A \to M_n(Q(B))$ with an
isomorphism $ M_n(Q(B)) \simeq Q(B)$ which is canonical in the sense that
it arises from an isomorphism $M_n(B) \simeq B$, we obtain a unital
extension $e(v) : A \to Q(B)$ of $A$ by $B$. By use of a unitary lift of
$\left( \begin{smallmatrix} v  & \\ & v^*\end{smallmatrix} \right)$ one
sees that $e(v) \oplus e(v^*)$ is split, proving that $e(v)$ represents an
element in $\Ext^{-1}_{\text{unital}}(A,B)$. If $v_t, t \in [0,1]$, is a
norm-continuous path of unitaries in $M_n(Q(B))$ there is a partition $0
=t_0 < t_1< t_2 < \dots < t_N =1$ of $[0,1]$ such that
$v_{t_i}v_{t_{i+1}}^*$ is in the connected component of $1$ in the unitary
group of $M_n(Q(B))$ and hence has a unitary lift to $M_n(M(B))$. It
follows that $e\left(v_0\right) = e\left(v_1\right)$, and it is then clear
that the construction gives us a group homomorphism
 $$
u : K_1\left( Q(B)\right) \to \Ext^{-1}_{\text{unital}}(A,B) .
 $$

\begin{lemma}\label{B35} The sequence
\begin{equation}\label{B36}
\begin{xymatrix}{
K_1(Q(B)) \ar[r]^-{u}  & \Ext^{-1}_{\text{\emph{unital}}}(A,B) \ar[r] &
\Ext^{-1}(A,B)  \ar[d]^-{i^*} \\ K_1\left(\alpha(A)'\cap Q(B)\right)
\ar[u]^-{i^*} &  &     \Ext^{-1}(\mathbb C,B)}
\end{xymatrix}
\end{equation}
is exact.
\begin{proof}
Exactness at $K_1(Q(B))$: If $v$ is a unitary in $M_n\left(\alpha(A)'\cap
Q(B)\right)$, the extension $\Ad v \circ \left(1_n \otimes \alpha\right) =
1_n \otimes \alpha$ (of $A$ by $M_n(B)$) is split, proving that $u \circ
i^* = 0$. To show that $\ker u \subseteq \im i^*$, let $v \in Q(B)$ be a
unitary such that $u[v] = 0$. Then $\Ad v \circ \alpha \oplus \alpha$ is
unitarily equivalent to $\alpha \oplus \alpha$, which means that there is
a unitary $S \in M(M_2(B))$ such that
\begin{equation}\label{B6}
\Ad \left(\left( \id_{M_2(\mathbb C)} \otimes q\right)(S)\left(
\begin{smallmatrix} v & {} \\ {} & 1 \end{smallmatrix} \right)\right)
\left(\begin{smallmatrix}  \alpha(a) & {} \\ {} & \alpha(a)
\end{smallmatrix} \right )    =  \left ( \begin{smallmatrix}\alpha(a) & {}
\\ {} & \alpha(a) \end{smallmatrix} \right )
\end{equation}
for all $a \in A$. Since the unitary group of $M(M_2(B))$ is
normconnected by \cite{M} or \cite{CH}, the unitary $ \left (
\begin{smallmatrix} v & {} \\ {} & 1 \end{smallmatrix} \right )$ is
homotopic to $\left( \id_{M_2(\mathbb C)} \otimes q\right)(S)\left(
\begin{smallmatrix} v & {} \\ {} & 1 \end{smallmatrix} \right)$ which is
in $M_2(\alpha(A)^{\prime} \cap Q(B))$ by (\ref{B6}). This implies that
$[v] \in \im i^*$. The same argument works when $v$ is a unitary
$M_n(Q(B))$ for some $n \geq 2$.

Exactness at $ \Ext^{-1}_{\text{unital}}(A,B)$: For any unitary $v \in
Q(B)$,
 $$
(\Ad v \circ  \alpha) \oplus 0 = \Ad \left( \id_{M_2(\mathbb C)}
\otimes q\right)(T) \circ ( \alpha \oplus 0),
 $$
where $T \in M_2(M(B))$ is
a unitary lift of $ \left ( \begin{smallmatrix} v & {} \\ {} & v^*
\end{smallmatrix} \right )$. Hence $[\Ad v \circ  \alpha] = 0$ in
$\Ext^{-1}(A,B)$. The same argument works when $v$ is a unitary
$M_n(Q(B))$ for some $n \geq 2$, and we conclude that the composition
$K_1(Q(B)) \to \Ext^{-1}_{\text{unital}}(A,B) \to  \Ext^{-1}(A,B)$ is
zero. Let $\varphi : A \to Q(B)$ be a unital extension such that
$[\varphi] = 0$ in $\Ext^{-1}(A,B)$. By Lemma \ref{B31} this means that
there is a unitary $T \in M(M_3(B))$ such that
 $$
\Ad \left( \id_{M_3(\mathbb C)} \otimes q\right)(T) \circ  \left (
\begin{smallmatrix} \varphi & {} & {} \\ {} &  \alpha & {} \\ {} & {} & 0
\end{smallmatrix} \right ) =  \left ( \begin{smallmatrix}  \alpha & {} &
{} \\ {} &  \alpha & {} \\ {} & {} & 0  \end{smallmatrix} \right ) .
 $$
It follows that $\left( \id_{M_3(\mathbb C)} \otimes q\right)(T) =  \left
( \begin{smallmatrix} V & {}  \\ {} & r  \end{smallmatrix} \right )$ for
some unitaries $V \in M_2(Q(B))$ and $r \in Q(B)$. Hence
 $$
\left(\begin{smallmatrix} \varphi & {}  \\ {} &  \alpha  \end{smallmatrix}
\right ) = \Ad V^* \circ \left ( \begin{smallmatrix}  \alpha & {}  \\ {}
& \alpha  \end{smallmatrix} \right ).
 $$
Thus $[\varphi] = u[V^*]$.

Exactness at $ \Ext^{-1}(A,B)$: It is obvious that $i^*$ kills the image
of $\Ext^{-1}_{\text{unital}}(A,B)$, so consider an invertible extension
$\varphi : A \to Q(B)$ such that $[\varphi \circ i] = 0$ in
$\Ext^{-1}(\Bbb C, B)$. By Lemma \ref{B31}, applied with $A =\mathbb C$,
this means that there is a unitary $T \in M_3(M(B))$ such that
\begin{equation}\label{B8}
\left(\id_{M_3(\mathbb C)} \otimes q\right)(T)\left ( \begin{smallmatrix}
\varphi(1) & {} & {} \\ {} & 1 & {} \\ {} & {} & 0  \end{smallmatrix}
\right )\left(\id_{M_3(\mathbb C)} \otimes q\right)(T^*) = \left (
\begin{smallmatrix} 0 & {} & {} \\ {} & 1 & {} \\ {} & {} & 0
\end{smallmatrix} \right ) .
\end{equation}
Set $\psi = \varphi \oplus
\alpha \oplus 0$. It follows from (\ref{B8}) that there are isometries
$W_1,W_2,W_3 \in M(B)$ and a unitary $u \in M(B)$ such that $W_i^*W_j = 0,
i\neq j$, $W_1W_1^* + W_2W_2^* + W_3W_3^* = 1$ and $\Ad q(u) \circ \psi(1)
= q(W_2W_2^*)$. Then $\Ad q(u) \circ \psi + \Ad q(W_1) \circ \alpha + \Ad
q(W_3) \circ \alpha$ is a unital extension which is invertible because it
admits a completely positive contractive lifting to $M(B)$ since $\psi$
does, cf. \cite{A}. As it represents the same class in $\Ext^{-1}(A,B)$ as
$\varphi$, the proof is complete.

\end{proof}
\end{lemma}

In order to complete the sequence of Lemma \ref{B35}, let $i_1^* :
\Ext^{-1}(A,B) \to K_1(B)$ be the composition
\begin{equation}\label{B61}
\begin{xymatrix}{
\Ext^{-1}(A,B) \ar[r] & K_1\left(\alpha(A)' \cap Q(B)\right)
\ar[r]^-{i^*}  & K_1\left(Q(SB)\right)) \ar[r] & K_1(B),}
\end{xymatrix}
\end{equation}
where the first map is the isomorphism (\ref{B32}) and the last is the
well-known isomorphism. Let $u_1 : K_1(B) \to
\Ext_{\text{unital}}^{-1}(A,SB)$ be the composition
\begin{equation}\label{B62}
\begin{xymatrix}{
K_1(B) \ar[r] & K_1\left( Q(SB)\right) \ar[r]^-{u}  &
\Ext_{\text{unital}}^{-1}(A,SB) ,}
\end{xymatrix}
\end{equation}
where the first map is the well-known isomorphism (the inverse of the one
used in (\ref{B61})) and second is the $u$-map as defined above, but with
$SB$ in place of $B$. Let $i_0^* : \Ext^{-1}(A,SB) \to K_0(B)$ be the
composition
\begin{equation}\label{B63}
\begin{xymatrix}{
\Ext^{-1}(A,SB) \ar[r]^{i^*} & \Ext^{-1}(\mathbb C, SB) \ar[r]  &
K_0\left(B\right),}
\end{xymatrix}
\end{equation}
where the second map is the well-known isomorphism. Finally, let $u_0
:K_0(B) \to \Ext_{\text{unital}}^{-1}(A,B)$ be the composition
\begin{equation}\label{B64}
\begin{xymatrix}{
K_0(B)   \ar[r] & K_1\left(Q(B)\right)  \ar[r]^-u  &
\Ext_{\text{unital}}^{-1}(A,B) ,}
\end{xymatrix}
\end{equation}
where the first map is the well-known isomorphism. We have now all the
ingredients to prove

\begin{thm}\label{cor1} The sequence
\begin{equation*}\label{diag6}
\begin{xymatrix}{
K_0(B) \ar[r]^-{u_0} &  \Ext^{-1}_{\text{\emph{unital}}}(A,B) \ar[r] &
\Ext^{-1}(A,B) \ar[d]^{i_1^*} \\ \Ext^{-1}(A,SB) \ar[u]^-{i_0^*}  &
\Ext^{-1}_{\text{\emph{unital}}}(A,SB) \ar[l] & K_1(B) \ar[l]^-{u_1} }
\end{xymatrix}
\end{equation*}
is exact.
\begin{proof} If we apply Lemma \ref{B35} with $B$ replaced by $SB$ we
find that the sequence
\begin{equation}\label{B65}
\begin{xymatrix}{
\Ext^{-1}(\mathbb C, SB)  &   &  K_1\left(\alpha(A)'\cap Q(SB)\right)
\ar[d]^-{i^*} \\ \Ext^{-1}(A,SB)  \ar[u]^-{i^*}  &
\Ext_{\text{unital}}^{-1}(A,SB) \ar[l] &  \ar[l]^-u K_1\left(Q(SB)\right)
}
\end{xymatrix}
\end{equation}
is exact. Thanks to the commuting diagram (\ref{B34}) we can patch this
sequence together with the sequence from Lemma \ref{B35} with the stated
result.
\end{proof}
\end{thm}

\end{document}